\documentstyle[12pt]{article}
\def\ba{\begin{array}}
\def\ea{\end{array}}	
\def\be{\begin{equation}}
\def\ee{\end{equation}}
\def\bem{\begin{em}}
\def\eem{\end{em}}
\def\ot{\otimes}

\def\ld{\ldots}

\def\a{\alpha}                           
\def\b{\beta}
\def\g{\gamma}
\def\G{\Gamma}
\def\d{\delta}
\def\D{\Delta}

\def\vp{\varphi}

\def\g{\gamma}
\def\G{\Gamma}
\def\si{\sigma}
\def\vp{\varphi}
\def\th{\Theta}
\def\om{\omega}
\def\Om{\Omega}
\def\ra{\longrightarrow}

\def\ca{{\cal A}}

\def\ca{{\cal A}}

\def\cw{{\cal W}}

\def\cm{{\cal M}}
\def\dn#1,{_{({#1})}}
\def\br{[.,.]_c}
\def\real{I\!\! R}

\def\1n{^{(1)}}
\def\2n{^{(2)}}
\def\3n{^{(3)}}
\begin{document}
\title{PARTICLES IN SINGULAR MAGNETIC FIELD\footnotemark}
\author{W{\l}adys{\l}aw Marcinek
\\ Institute of Theoretical Physics, University
of Wroc{\l}aw,\\ Pl. Maksa Borna 9, 50-204  Wroc{\l}aw,\\
Poland}
\date{}
\maketitle
\footnotetext[1]
{The paper is partially supported by KBN, Grant No 2P 302 087 06\\
and by the Katholischer Akademischer Auslander Dienst, Bonn, Germany\\
e-mail: wmar@ift.uni.wroc.pl}
\begin{abstract}
\noindent
An algebraic formalism for description of quantum states of charged 
particle with spin moving in two--dimensional space under influence of 
singular magnetic field is developed in terms of graded algebras. The 
fundamental assumption is that the particle is transformed into a composite 
system which consists quasiparticles, quasiholes and magnetic fluxes. Such 
system is endowed with generalized statistics determined by a grading 
group and a commutation factor on it. Composite systems corresponding 
to the quantum Hall effect and the electronic magnetotransport anomaly 
are described. The Fock space representation are also given.
\end{abstract}

\vspace{3.5cm}

PACS number classification: 03.65.Fd, 12.35.K, 11.30.Pb, 02.40+m\\

\newpage
\section{Introduction}
In the last years charged particles moving under influence of 
a strong and singular magnetic field have been studied from a 
few different points of view by many authors. New and very 
specific effects in two dimensional systems of particles
have been discovered, 
see Ref. \cite{tsg,eha,hal,jai,dst,ara} for example. 
A fundamental quantity for the description of these new effects
is the so--called {\it filling factor} $v$ which is defined as the
ratio $\frac{p}{N}$ of the number of particles $p$ by the number
of magnetic fluxes $N$. {\it A magnetic flux} is a magnetic field 
concentrated completely on a vertical line carrying an elementary
flux $\Phi_0 := \frac{h}{e}$. For integer--valued $v$ there is the
well--known {\it integer quantum Hall effect} (IQHE). For 
fractional--valued $v$ with odd denominator there is the famous 
{\it fractional quantum Hall effect} (FQHE) \cite{tsg}. 
It is very interesting and
fascinating fact that for fractional--valued filling factor with 
even denominator the corresponding effect is completely different 
from the mentioned above quantum Hall effect. This is the so--called
{\it electronic magnetotransport anomaly} (the $\frac{1}{2}$--anomaly).
In this case the magnetic field is completely compensated such
that particles move like in the the absence of the magnetic field!
In this way we obtain the following scheme
\be
v = 
\left\{
\ba{ll}
1, 2, \ldots&\mbox{IQHE},\\
\\
\frac{1}{2}, \frac{1}{4}, \ldots &\mbox{the} \; \frac{1}{2}\; 
\mbox{anomaly},\\
\\
\frac{1}{3}, \frac{1}{5}, \ldots&\mbox{FQHE}.
\ea
\right.
\ee
At present there is no satisfactory theory for unified
description of all these effects.
The most popular theory for the quantum Hall effect is the
Chern\--Simons\--Landau\--Ginzburg (CSLG) theory, see for
example \cite{kaz} and reference therein. In this theory
electrons are described as bosons carrying odd number of flux
quanta. Such composite system is said to be {\it a composite boson}.
The $\frac{1}{2}$--anomaly can be explained in a similar way.
In this case electrons can be considered as a degenerate system 
of new fermions in the absence of magnetic field. Such particles 
are called {\it composite fermions} \cite{jai,dst,gsj}. They can be 
imagined as electrons carrying even number of magnetic fluxes. 

An interesting concept for composite system which consist a particle 
with charge $e$ and a singular magnetic field concentrated completely 
on a vertical line has been introduced by Wilczek \cite{wil}, see 
also \cite{wiz}. The charge is moving around 
the singular line in the flat region in which there is no influence 
of the magnetic field. Every rotation yields certain phase factor $q$.
Observe that the space on which the particle moving is 
$\cm \equiv I\!\! R^2 \setminus \{s \}$, where $s \in I\!\! R^2$ 
is the point of intersection of the magnetic line with the plane 
$I\!\! R^2$. The factor $q$ can be in general a complex number of 
the following form: $q := \exp(i \vp)$, where $0 \underline{<} \vp < 2\pi$ 
is the so--called statistics parameter. This means that the
statistics of the particle is determined by the value of $\vp$.
For $\vp = 0$ we have boson, for $\vp = \pm \pi$ -- fermion.
For arbitrary $\vp \in [0, 2\pi)$ we obtain anyons. Note that 
the name "anyons" is sometimes reserved to the case when $q$ 
is the $m$-th root of unity, see for example the papers of Majid 
\cite{anm,bnm}, where anyons are given in terms of certain 
$Z_m$-graded structures. The description of the statistics 
corresponding to composite fermions or bosons is a problem. 
A few remarks on this problem has been given by Jacak, Sitko 
and Wieczorek in Ref.\cite{jsw}. According to their suggestion 
the value of the statistics parameter not need to be restricted
to the interval $[0, 2\pi)$, but it can be an arbitrary real number. 
If $\vp = (2k+1)\pi$, $k=1, 2,\ldots$ for example, then one 
can expect that we obtain nothing else but the usual fermion statistics. 
But according to  Jacak, Sitko and Wieczorek in this case we obtain new 
kind of statistics, just the statistics of composite fermions. The case 
of $\vp = 2k \pi$ $k=1, 2,\ldots$ corresponds to composite bosons. 
An algebraic formalism for particle system with generalized 
statistics has been considered by the author in Ref. \cite{WM4,WM6,wm7}
see also \cite{wmq,qweyl,top,mco}. It is interesting that in this
algebraic approach all commutation relations for particles equipped 
with arbitrary statistics can be described as a representation of the 
so-called quantum Weyl algebra $\cw$ (or Wick algebra)
\cite{RM,m10,jws,ral}. Similar approach has been also considered by
others authors, see \cite{sci,twy,mep,jws} and \cite{mphi,mira}.
In this attempt the creation and
anihilation operators act on certain quadratic algebra $\ca$. 
The creation operators act as the multiplication in $\ca$ and the 
anihilation ones act as a noncommutative contraction (noncommutative
partial derivatives). The algebra $\ca$ play the role of noncommutative 
Fock space. The applications for particles
in singular magnetic field has been described in Ref. \cite{top,mco}.
As a result the fact that at $v = \frac{1}{2}$ the system
is degenerate and at $v = \frac{1}{N}$ for $N$ odd the Landau level
is fractionally filled has been obtained in an obvious way. 

In this paper an algebraic formalism for a system of charged 
particles with spin $\frac{1}{2}$ in singular magnetic field
is developed. We restrict our attention to the system with the 
filling factor $v = \frac{1}{N}$ for simplicity. 
The generalization for arbitrary filling factor is possible 
but it need more complicated notation. 
The fundamental assumption of the paper is 
that the particle is transformed into a composite system which 
consists quasiparticles, quasiholes and magnetic fluxes. 
The transformation a result of interactions of charged particle 
with spin with magnetic field.
If $N$ is even, then the system is identified as composite 
fermions and if $N$ is odd, then we obtain composite bosons.
It is obvious the statistics of our composite systems can not 
be described by the scalar statistics parameter $q = \exp (i \vp)$.
We use here the concept of commutation factor \cite{sch}
in order to describe the statistics.
The paper is organized as follows.
The fundamental concept of the paper based 
on Refers \cite{qweyl,top,mco} is given In Section 2.
In Section 3 the quantum Weyl algebra 
is described in details. It is shown that it is in fact
a $\G$--graded generalized (color) Lie algebra, 
\cite{sch}. The representation of the quantum Weyl algebra
is considered in Section 4. The Fock space of quantum states
for composite fermions and composite bosons is described
in Section 5. The notion of commutation factors is shortly 
summarized in the Appendix.
Note that the algebraic approach to flat particle system in singular 
magnetic field is not complete theory but it can be further developed.
One can write down the Hamiltonian and study the corresponding field
theory. One can use the formalism described by Liguori and Mintchev
\cite{lim} for such study .
We hope that the possibility for the unified description 
of the quantum Hall effect and the electronic magnetotransport 
anomaly based on the presented in this paper algebraic formalism
should be clear.
\section{Fundamental assumptions}
Let us consider a system of charged particles with spin $\frac{1}{2}$ 
moving on a flat space $I\!\! R^2$ under influence of singular 
perpendicular magnetic field. The magnetic field is sufficiently strong 
for the polarization effects of the spin. The {\it singularity} of the 
magnetic field means here that the field is completely concentrated 
into fluxes. We assume that there is in average $N$ magnetic fluxes 
per particle. This means that the filling factor is $v = \frac{1}{N}$.
Our fundamental assumption is that the problem of many particles
interacting with the singular magnetic field can be 
reduced to the study of a system consisting just one  
particle and $N$ magnetic fluxes. In this way we can 
restrict our attention to study of such system. 
It is known that quantum states of charged particle 
in the magnetic field are described by Landau levels. 
Here is an additional degeneracy of such states 
connected with the singular nature of the field.
The structure of such degeneration is determined
by the change of phase of the particle.
Let us study the degeneration in more details.
We assume that the particle is in the Landau lowest level.
The "effective" space for the particle in singular magnetic field is
$\cm \equiv \real ^2 \setminus \{s_1,...,s_N\}$, where $s_i$ is
the point of intersection of the $i$-th magnetic flux with the
plane. Let us denote by $P_m\cm$ the space of all homotopy classes of 
paths which starts at $m_0$ and and end at arbitrary point $m \in \cm$. 
It is known that the union $\cup_m P_m\cm = P\cm$ is a covering space 
$P = (P\cm, \pi, \cm)$. The projection $\pi : P\cm \ra \cm$
is given by
\be
\ba{ccc}
\pi (\xi) = m &\mbox{iff}&\; \xi \in P_m \cm .
\ea
\ee
The homotopy class $P_{m_0}\cm$ of all paths which start at $m_0$
and end at the same point $m_0$, (i.e. a loop space) can be
naturally endowed with a group structure. This group is known
as the fundamental group of $\cm$ at $m_0$ at $m_0$ and is denoted 
by $\pi_1 (\cm, m_0)$. Observe that two groups $\pi_1 (\cm, m_0)$
and $\pi_1 (\cm, m)$ at two arbitrary points $m_0$ and $m$, respectively,
are isomorphic. Hence we can introduce the fundamental group 
$\pi_1 (\cm)$ for the whole space $\cm$ in an obvious way.
Note that in our case the fundamental group $\pi_1 (\cm)$ 
for the space $\cm \equiv \real ^2 \setminus \{s_1,...,s_N\}$ 
is a free group generated by $\tilde{\si}_i$, $(i = 1,...,N)$ 
and is denoted by 
$G_N \equiv G_N (\tilde{\si}_1 ,\ldots, \tilde{\si}_N )$, 
where $\tilde{\si}_i$ is the homotopy class of all paths which 
encloses the point $m_i$ and none of the remaining points
$m_j$ for $j \neq i$, see Ref.\cite{tqst} for example for more details.
Next we assume that for every class of loops the phase 
change of particle under consideration is fixed.
We introduce an equivalence relation in the fundamental group
$G_N$. Two homotopy classes of loops in $G_N$ are equivalent if 
and only if the phase changes corresponding to this classes are equal.
The equivalence class of $\tilde{\si}_i$ with respect to the above
relation is denoted by $\si_i$. Let $\G$ be a quotient group of the 
fundamental group $G_N$ by the above relation. It is obvious that 
the phase change corresponding for $\tilde{\si}_i \tilde{\si}_j$
and for $\tilde{\si}_j \tilde{\si}_i$ is the same for all $i \neq j$.
Hence we have the relation $\si_i \si_j = \si_j \si_i$ in the quotient
group $\G$. Observe that generators $\si_i$ of $\G$ can be identified
with rotations around axis at the point $s_i$. This means that the 
group $\G$ can be identified with the group $Z \oplus\ldots\oplus Z$ 
which contains $N$ copies of the group of integers. Every rotation
around the axis at $s_i$ yields the same phase factor $q_i$. 
Two rotations around two different axis $s_i$ and $s_j$ are
independent and they yields the phase factor $c_{ij}$. It is
natural to assume that $c_{ij} c_{ji} = 1$ and $c_{ii} \equiv q_i$.
It is obvious that there is a commutation factor 
$c : \G \times \G \ra I\!\!\!\! {\bf C} \setminus \{0\}$
on the group $\G$ (see the Appendix) such that 
$c(\si_i, \si_j) \equiv c_{ij}$. The commutation factor $c$ 
generalize the scalar statistics parameter $q = \exp (i \vp)$, where 
$\vp$ is an arbitrary real number. 
It follows from the above considerations that
the structure of degenerations of Landau levels can be described 
by the commutation factor $c$ on $\G$.
Structures of different kind of degenerations
correspond to nonequivalent 
classes of commutation factors on the group $\G$. 
This means that we can use the concept of color statistics and 
related $\G$--graded structures in order to obtain the algebraic 
formalism for the particle in the singular magnetic field \cite{wm7}.
In is interesting that in certain particular cases, see the Appendix,
the grading group can be reduced to the group 
$\G = Z^N_n \equiv Z_n \oplus...\oplus Z_n$ for $n >2$ or to
$Z^N_2 \equiv Z_2 \oplus...\oplus Z_2$. 

It is interesting that the structure of degeneration of the 
Landau levels can be described as quantum states of new
particles equipped with generalized statistics. 
We assume that the particle under considerations is transformed 
into a composite system which consists quasiparticles, quasiholes
and magnetic fluxes. Such composite system is said to be 
{\it a composite particle or generalized quon}. 
The transformation of the particle into the composite system is 
a result of interactions of charged particles with spin with 
magnetic field. Note that if the particle is coupled with certain 
flux in such a way that the magnetic field of this flux is canceled, 
then  we say that we have {\it a quasiparticle}. 
In this case the particle is said to be {\it bound} to the flux.
In the opposite case, i.e. when the magnetic field is not canceled 
we say on quasiholes. {\it A quasihole} is a "free" flux which behave 
like particle endowed with fermion statistics. Hence the points of 
intersections of $N$ magnetic fluxes with the plane must be a set
of $N$ different points $s_1,...,s_N$.
The number of quasiparticles and quasiholes is 
equal to the number of magnetic fluxes. This means that the filling
factor for the transformed system is $v' = \frac{m+n}{N} \equiv 1$,
where $m$ is the number of quasiparticles, and $n$ is the number
of quasiholes. The "effective" magnetic field is
\be
B_{eff} := B - m \Phi_0 ,
\ee
where $B$ is the external magnetic field. Observe that if 
$B = N \Phi_0$ and the number of quasiparticles is equal to $N$,
then the magnetic field is completely canceled!
A quasiparticle is in fact the charged particle bound to single
magnetic flux. The particle bound to two different magnetic fluxes
are understand as a system of two different quasiparticles. Obviously 
the particle can not be bound to two fluxes at the same point! 
Quasiparticles as components of certain composite particle 
have also their own statistics. This statistics is determined
by the commutation factor $c$ on the group $\G$.
For electron in singular magnetic field it is natural to assume 
that we have $c_{ii} \equiv \exp (i\vp)$, where $\vp = \pi (N+1)$, 
i. e. $c_{ii }\equiv -(-1)^N$ and 
$c_{ij} \equiv \exp (i\pi N) \equiv (-1)^N$ for $i \neq j$.
The reason that in the formula for $c_{ii}$ we have that
$\vp$ is equal to the number of magnetic fluxes $N$ plus $1$
is the spin of electrons. The above expressions for the factor 
$c$ means that two different quasiparticles for $N$ even commute 
and for $N$ odd -- anticommute. This also means that the composite 
boson can contain only one quasiparticle. The number of quasiparticles
can be equal to $N$ for composite fermions only! In others words,
our quasiparticles are examples of particles equipped with the
so--called color statistics, i.e. the parastatistics in Green
representation \cite{wm7}. This means that one can use the
concept of color statistics and related structures with
abelian grading group in order to obtain an algebraic 
formalism for particles in singular magnetic field.	
\section{On quantum Weyl algebras}
Let us denote by $\th^i$ the quantum quasiparticle state 
coupled with the flux at $s_i$. The corresponding 
conjugate state is denoted by $\th^{\ast}_j$.
It is natural to assume that there is a finite dimensional vector 
space $E$ equipped with a basis $\th^i, i = 1,...,N = dim E$. The complex 
conjugate  space $E^*$ is endowed with a basis $\th^{\ast}_{i}$ such that
we have that $<\th_i^{\ast}|\th^j> = \d_i^j$ for $i, j = 1, \ldots , N$.
Let us give the notion of quantum Weyl algebras 
$\cw \equiv \cw_{\G, c}(N)$ corresponding to an arbitrary
grading group $\G$ and a commutation factor $c$ on it.\\
\begin{em}
{\bf Definition:}
A $(\G,c)$-quantum Weyl algebra is a quotient algebra 
\be
\cw \equiv \cw_{\G,c}(N) := T(E \oplus E^*)/I_{\G,c},
\ee
where $I_{\G,c}$ is an $*$-ideal in the tensor algebra $T(E \oplus E^*)$
generated by the following elements 
\be
\ba{l}
\th^*_i \ot \th^j - c_{ij} \ \th^{j} \ot \th^*_i - \d^j_i {\bf 1},
\;\;\;\mbox{for all}\;\; i, j,\\
\th^i \ot \th^j - c_{ij} \ \th^j \ot \th^i, \;\;\;
\th^*_j \ot \th^*_i - c_{ij} \ \th^*_i \ot \th^*_j,
\;\;\;\mbox{for}\; i \neq j,\\
(\th^i)^2, \;\;\; (\th_j^{\ast})^2 \;\;\; \mbox{if}\;\;c_{ii} = -1,
\ea
\ee
(no sum), for $i, j = 1,...,N$, where $c_{ij}$ are coefficient 
of certain commutation factor on the grading group $\G$.\\
\end{em}
The above definition means that $\cw_{\G, c}(N)$ is a $*$-algebra 
generated by $\th^i$ and $\th^{\ast}_j$; $(i,j = 1,...,N)$ subject 
to relations
\be
\ba{llll}
\th^*_i \th^i = 1 + q_i \ \th^i \th^*_i&\;\;\mbox{for}\;\;i = j,&
\th^*_i \th^j = c_{ji} \ \th^j \th^*_i&\;\;\mbox{for}\;\;i \neq j,\\
\th^i \th^j = c_{ij} \ \th^j \th^i&\;\;\mbox{for all}\;\; i,j,&
\th^*_j \th^*_i = c_{ij} \ \th^*_i \th^*_j&\;\;
\mbox{for all}\;\;i,j,\\
(\th^i)^2 = 0,&(\th_j^{\ast})^2 = 0,& \mbox{if}\;c_{ii} = -1,&
\label{crul}
\ea
\ee
where $q_i := c_{ii} : i = 1,\ldots, N$ are diagonal elements of $c$,
note that the same symbols for elements in the tensor algebra 
$T(E \oplus E^*)$ and for generators of the algebra $\cw_{\G, c}(N)$
have been used for simplicity. The algebra $\cw_{\G, c}(N)$ should
be also denoted by 
$\cw_{\G, c}<\th^1,\ldots,\th^N,\th^*_1,\ldots,\th^*_N >$.\\
\begin{em}
{\bf Theorem:} The quantum Weyl algebra $\cw \equiv \cw(\G, c)$ 
is a $\G$--graded $c$--Lie algebra.\\
\end{em}
{\bf Proof:} 
We introduce the $\G$--gradation of the algebra $\cw$ as follows:
for generators of the algebra we define that
$\mbox{grade}(\th^i) \equiv |\th^i| := \si_i$, and
$\mbox{grade}(\th_i^*) \equiv |\th_i^*| := -\si_i$,
where $\{\si_i\}_{i=1}^N$ is a set of generators of $\G$.
All monomials in $\th^i$  and $\th^{\ast}_j$ modulo
generating relations are homogeneous elements of $\cw$.
We use the formula $|XY| = |X| + |Y|$ for the extension of
gradation for arbitrary homogeneous elements of the algebra.
The generalized bracket $\br : \cw \ot \cw \ra \cw$
is defined by the formula
\be
\br X \ot Y \equiv [X, Y]_c := X Y - c(|X|, |Y|) Y X,
\ee
where $X, Y$ are arbitrary homogeneous elements of the algebra $\cw$.
We can calculate the bracket as follows: for generators we obtain
the relation
\be
\ba{l}
[\th_i^{\ast}, \th^j]_c
:= \th_i^{\ast} \th_j - c(|\th_i^{\ast}|, |\th^j|) \th^j \th_i^{\ast}
= \d_i^j ,
\ea
\ee
where $c(|\th_i^{\ast}|, |\th^j|) = c(-\si_i, \si_j) = c_{ji}$
and the relations (\ref{crul}) have been used. In the similar way
we obtain the formula
\be
\ba{l}
[\th^i, \th_J^{\ast}]_c
:= \th^i \th^{\ast}_j - c(|\th^i|, |\th^{\ast}_j|) \th^{\ast}_j 
\th^i\\ 
= - c(|\th^i|, |\th^{\ast}_j|) 
(\th_i^{\ast} \th_j - c(|\th_i^{\ast}|, |\th^j|) \th^j \th_i^{\ast})
= - c_{ij} \d_i^j .
\ea
\ee
For the extension of our generalized bracket to the higher order
homogeneous elements of the algebra we can use the following formulae
\be
\ba{l}
[\th^{\ast}_i , \th^j X]_c =  
c_{ij} \th^j [\th^{\ast}_i, X] + \d_i^j X,\\
\left[
\theta_i^{\ast} X, \th^{j}]_{c} = \th^{\ast}_i [X, \th^{j}
\right] 
 - c_{ji} c(|X|, \si_j ) \d_i^j X,
\ea
\ee
for arbitrary homogeneous $X$ in $\cw$.
One can prove these formulae by induction. One can see that
the generalized bracket is $c$--anticommutative
\be
\ba{l}
[X, Y]_c = - c(|X|, |Y|) [Y, X]_c
\ea
\ee
for homogeneous $X$ and $Y$ in $\cw$. Finally, it not difficult
to calculate the following generalized Jacobi identity
\be
\ba{l}
[X, [Y, Z]_c ]_c = [[X, Y]_c , Z]_c
 + c(|X|, |Y|) \ [Y, [X, Z]_c ]_c .
\ea
\ee
\hfill $\Box$\\
It is interesting that for arbitrary generalized Lie algebra 
$L \equiv L_{\G, c}$ there is a corresponding Lie algebra or
superalgebra $s(L)$, \cite{sch,WM3}. 
This algebra is called the superisation 
of $L$ and we have the following crossed product
\be
\ba{l}
L := s(L) =\!\!\!\!\!\!/\!\!/ \ C_{\G, b}(N),
\label{cpo}
\ea
\ee
where $C_{\G, b}(N)$ is $\G$--graded $b$--commutative $\ast$--bialgebra.
This means that we have the relation $gh = b(|g|, |h|)hg$ for all 
homogeneous $g, h \in C_{\G, b}(N)$ of grade $|g|$ and $|h|$, 
respectively. The factor $b$ is given by the formula (\ref{ncf}). 
The bialgebra $C_{\G, b}(N)$ as an unital, associative algebra 
is generated by $e^i , e^{\ast}_j$ and the following relations
\be
\ba{ll}
e^i e^j = b_{ij} e^j e^i ,\;\;
e^{\ast}_i e^{\ast}_j = b_{ij} e^{\ast}_j e^{\ast}_i ,&
\;\;\mbox{for}\;\;i \neq j,\\
e^{\ast}_i e^j = e^j e^{\ast}_i = \d_i^j,&\;\;\mbox{for all}\;i, j,\\
\ea
\ee
(no sum). The gradation is given
in the standard way. This means that $|e^i| = \si_i$ and 
$|e^{\ast}_i| = - \si_i$, $\si_i$ are generators of the group $\G$.
The $\ast$--operation and the comultiplication by the formulae
\be
(e^i)^{\ast} := e^{\ast}_i, \;\;\;
(e^{\ast}_i)^{\ast} := e^i, 
\ee
i. e. $(e^i)^{\ast} \equiv (e^i)^{-1}$ and
\be
\D (e^i) := e^i \ot 1 + 1 \ot e^i,
\ee
respectively. The extension to the whole bialgebra is obvious \cite{mkl}.
For the tensor product we have here the relation
\be
(g \ot h)(k \ot l) = b(|h|, |k|) gk \ot hl
\ee
for homogeneous $g, h, k, l \in C_{\G, b}(N)$. Observe that the
bialgebra $C_{\G, b} (N)$ is not an usual Hopf algebra, but
a $\G$--graded generalized Hopf algebra, see \cite{wma}.
The relation (\ref{cpo}) means that $L$ is a subalgebra
of the tensor product $s(L) \ot C_{\G, b}(N)$ spanned by elements
of the form $X = x \ot g$ for homogeneous $X \in L$, $x \in s(L)$
$g \in C_{\G, b}(N)$ of the same grade, i.e. 
$|X| = |x| = |g| = \a \in \G$.
The algebra $s(L)$ is also $\G$--graded, but this gradation can
be reduced to the group $\G/\G_0$. The reduction is given by the
quotient map $\pi : \G \ra \G/\G_0)$. In this way we obtain that
the gradation of the algebra $s(L)$ can be given by the group
$Z_2$ or is trivial. Let us describe the superisation 
$s(\cw_{\G, c}(N))$ of the quantum Weyl algebra 
$\cw_{\G, c}(N) \equiv \cw_{\G, c}<\th^1 ,\ldots, \th^N, 
\th^{\ast}_1,\ldots, \th^{\ast}_N >$ in terms of generators. 
The algebra $s(\cw_{\G, c}(N))$ is generated by $x^i$ and $x^{\ast}_j$
$(i, j = 1, \ldots , N)$ and relations
\be
\ba{ll}
x^{*}_i x_i = 1 + q_i \ x_{i} x^{\ast}_i&
\;\;\mbox{for}\;\;i = j,\\
x^{\ast}_i x_j = c'_{ji} \ x_j x^{\ast}_i&\;\;\mbox{for}\;\;i \neq j,\\
x^i x^j = c'_{ij} \ x^j x^i&\;\;\mbox{for all}\;\; i,j.
\ea
\ee
where $c'$ is given by the formula (\ref{sad}).
This means that $s(\cw_{\G, c}(N))$ is also a quantum Weyl 
algebra, $s(\cw_{\G, c}(N)) \equiv \cw_{\pi(\G), c'}$.
The crossed product of algebras $s(\cw_{\G,c}(N))$ 
and $C_{\G, b}(N)$ is a subalgebra of the tensor product 
$s(\cw_{\G, c}(N)) \ot C_{\G, b}(N)$ generated by elements
\be
\ba{ccc}
\th^i := x^i \ot e^i&\th_i^* := x_i^{\ast} \ot e^{\ast}_i&
\;\mbox{(no sum)},
\label{sge}
\ea
\ee
Note that the above expressions for $\th^i$ and $\th^{\ast}_j$
are not unique. There is a freedom of choosing of the generators
of the bialgebra $C_{\G, b}(N)$ described by
the orthogonal group $O(N)$, \cite{wme}.
Let us consider a few simple examples:\\
{\bf Example 1}
Assume that $\G = Z \oplus \ldots \oplus Z$ ($N$--sumands) and
the factor $c$ is given by the relation (\ref{omo}), where
$\Omega_{ij} = 1$ for $i \neq j$.
Hence the commutation rules (\ref{crul}) can be given in the following 
form
\be
\ba{ll}
\th^*_i \th^i - q_i \th^i \th^*_i = 1&\;\;\mbox{for}\;\;i = j,\\
\th^*_i \th^j - \omega^{-1} \ \th^j \th^*_i = 0&\;\;\mbox{for}\;\;i < j,\\
\th^i \th^j - \omega \ \th^j \th^i = 0&\;\;\mbox{for}\;\; i < j.
\ea
\ee
{\bf Example 2} Let us consider the case corresponding for 
$\G = Z_2 \oplus \ldots \oplus Z_2$, in more detail. 
If we substitute the commutation factor $c$ of the form (\ref{com})
into the formulae (\ref{crul}), then we obtain the following relations
\be
\ba{ll}
\th^*_i \th^i +(-1)^N \th^i \th^*_i = 1&\;\;\mbox{for}\;\;i = j,\\
\th^*_i \th^j - (-1)^{\Om_{ij}} \ \th^j \th^*_i = 0&
\;\;\mbox{for}\;\;i < j,\\
\th^i \th^j - (-1)^{\Om_{ij}} \ \th^j \th^i = 0&\;\;\mbox{for}\;\; i < j.
\ea
\ee
It is interesting that in this case the bialgebra $C_{\G, b}(N)$
as algebra is reduced to the standard Clifford algebra $C_N$. In fact
we have $(e^i)^{\ast} \equiv e^i$ and we obtain the well-known
relations $e_i e_j + e_j e_i = 2 \d_{ij}$ for $i, j = 1, \ldots, N$.\\
{\bf Example 3}
If $N$ is even, we obtain the following relations for composite fermions
\be
\ba{ll}
\th^*_i \th^i - \th^i \th^*_i = 1&\;\;\mbox{for}\;\;i = j,\\
\th^*_i \th^j + \th^j \th^*_i = 0&\;\;\mbox{for}\;\;i < j,\\
\th^i \th^j + \th^j \th^i = 0&\;\;\mbox{for}\;\; i < j.
\ea
\ee
The quantum Weyl algebra generated by these relations is denoted by 
$\cw_{cf}(N)$. We can see that the algebra $s(\cw_{cf}(N))$ is
generated by the usual canonical commutation relations for the
system of $N$ bosons.\\
{\bf Example 4}
If $N$ is odd, then we obtain relations for composite bosons
\be
\ba{ll}
\th^*_i \th^i + \th^i \th^*_i = 1&\;\;\mbox{for}\;\;i = j,\\
\th^*_i \th^j - \th^j \th^*_i = 0&\;\;\mbox{for}\;\;i < j,\\
\th^i \th^j - \th^j \th^i = 0&\;\;\mbox{for}\;\; i < j.
\ea
\ee
The corresponding quantum Weyl algebra is denoted by $\cw_{cb}(N)$. 
In this case the algebra $s(\cw_{cb}(N))$ is generated by the usual 
canonical anticommutation relations for the system of $N$ fermions.\\
{\bf Example 5} Let $C_N$ be the Clifford algebra.
If we substitute
$\th^i \equiv \th^{\ast} \equiv \frac{1}{2} e_i$, and
$$
c_{ij} = \left\{
\ba{lll}
+1&\mbox{for}&i = j,\\
-1&\mbox{for}&i \neq j ,
\ea
\right.
$$
then we obtain a particular example of quantum Weyl algebra
denoted by $C_{\G}(N)$.
\section{Representation of the quantum Weyl algebra}
In this section we are going to construct a representation of the 
quantum Weyl algebra $\cw \equiv \cw_{\G, c}(N)$ on the $c$--symmetric 
algebra $\ca$. 
An algebra defined as the quotient $\ca \equiv \ca_c (E) := TE/I_{c}$
where $I_{c}$ is an ideal in $TE$ generated by the following relations
\be
\ba{ll}
\th^i \ \th^j - c_{ij} \ \th^j \ \th^i ,&\mbox{for}\;\; i \neq j,\\
(\th^i)^2 &\mbox{if}\;\; c_{ii} = -1,
\label{ide}
\ea
\ee
is said to be $c$-symmetric algebra over $E$ \cite{WM3}.
The algebra $\ca \equiv \ca_c (E)$ can be also denoted by 
$\ca_c <\th^1,\ld, \th^N>$, sometime.
It is easy to see that $\ca_c (E)$ is
$\G$-graded $c$-commutative algebra\\
\bem
{\bf Theorem:} 
Let $\cw$ be a Clifford--Weyl algebra, $\ca$ be a $c$-symmetric
algebra, then there is a representation $a : \cw \ra End(\ca)$
of the quantum Weyl algebra $\cw \equiv \cw_{\G, c}(N)$ on the
$\G$--graded $c$--commutative algebra $\ca$.\\
\eem
{\bf Proof:} We must prove that there is a well--defined homogeneous 
homomorphism $a : \cw \ra End(\ca)$ which transform the algebra $\cw$ 
as generalized Lie algebra into the generalized Lie algebra $End(\ca)$ 
of linear endomorphisms of $\ca$. This means that we have the relation 
\be
a_{[X, Y]_{c}} f = [a_{X}, a_{Y}]_{c} f
\ee
for every homogeneous $X, Y \in \cw$ and $f \in \ca$. The bracket
in the right hand side of the above formula is defined by the 
relation
\be
[a_{X}, a_{Y}]_{c} := a_{X} a_{Y} - c(|X|, |Y|) a_{Y} a_{X}.
\ee
For generators $\th^i$ and $\th^{\ast}_j$, $(i, j = 1,...,N)$ we define
\be
\ba{c}
a_{\th^i} f := m(\th^i \ot f) \equiv \th^i f,
\;\;a_{\th^{\ast}_i} f \equiv ev_k (\th^{\ast}_i \ot f)
\label{cao}
\ea
\ee
for every $f \in \ca^k$, where 
$ev_k : E^* \ot \ca^{\ot k} \ra \ca^{\ot k-1}$ are a set of linear 
mappings befined by the following formulae
\be
\ba{l}
ev_1 (\th^{\ast}_i \ot \th^j) := <\th^{\ast}_i |\th^j > = \delta_i^j,\\
ev_{k} (\th^{\ast}_i \ot \th^j f) := ev_1 (\th^{\ast}_i \ot \th^j) f + 
c_{ij} \ \th^j ev_{k-1}(\th^{\ast}_i \ot f),
\label{eva}
\ea
\ee
for $f \in \ca^{k}$. These mappings are said to be 
{\it a right evaluation}. Observe that we have the 
relation
\be
\ba{c}
ev_{k-1} (\th^{\ast}_i \ot ev_k (\th^{\ast}_j \ot f))
- c_{ij} ev_{k-1} (\th^{\ast}_j \ot (ev_k (\th^{\ast}_i \ot f)) = 0.
\label{rco}
\ea
\ee
on $E^* \ot E^* \ot \ca^k$. We have for example
\be
ev_2 (\th^{\ast}_i \ot \th^k \th^l) = \delta_{i}^k \ \th^l
+ c_{ik} \ \delta_{il} \ \th^k .
\ee
We can see that the evaluation $ev := \{ev_k : k = 1, \ldots \}$ 
can be well defined in a consistent way on the whole algebra $\ca$, 
see \cite{wmq,RM,m10}. We can also see that the definition
(\ref{cao}) can be extended in a consistent way for arbitrary 
homogeneous element of the algebra $\cw$. Obseve that we have
here the following lemma:\\
\hfill $\Box$\\
\begin{em}
{\bf Lemma:}
We have on the algebra $\ca$ the following commutation
relations for the representation $a$ of algebra $\cw_{\G, c}(N)$
on $\ca$
\be
\ba{l}
[a_{\th^{\ast}_i}, a_{\th^j}]_c = \d_i^j {\bf 1},\;\;
\left[a_{\th^i}, a_{\th^j}\right]_c = 0,\;\;
[a_{\th^{\ast}_i}, a_{\th^{\ast}_j}]_c = 0.
\label{cr2}
\ea
\ee
\end{em}
{\bf Proof} Using the relations (\ref{cao}) for the 
first relation (\ref{cr2}) we obtain
\be
\ba{l}
\left(
[a_{\th^{\ast}_i}, a_{\th^j}]_c
\right) f =
\left(
a_{\th^{\ast}_i} \ a_{\th^j} - c_{ij} \ a_{\th^{\ast}_j} \ a_{\th^i}
\right) f\\ 
= \left[
ev_{l+1}(\th^{\ast}_i \ot \th^j  f) - c_{ij} \ 
 \th^j ev_l (\th^{\ast}_i \ot f)
\right]
= ev_1 (\th^{\ast}_i \ot \th^j) f = \delta_{ij} f ,
\ea
\ee
where $f \in \ca^l$ and the relation (\ref{eva}) has been used. 
The second relation (\ref{cr2}) follows immediately from the
$c$--commutativity of the algebra $\ca$. The last relation follows
from the equation (\ref{rco}).
\hfill $\Box$\\
Let $a : \cw_{\G, c}(N) \ra End(\ca)$ be a Fock representation
of $\cw_{\G, c}(N)$ on $\ca$. 
Then there is a corresponding Fock representation
$\tilde{a} : s(\cw_{\G, c}(N)) \ra End(s(\ca_c))$
of the algebra $s(\cw_{\G, c}(N))$ on $s(\ca)$, where
$s(\ca)$ is a commutative (or supercommutative) algebra,
the superisation of $\ca$.
For the representation $a$ we have the following decomposition
\be
\ba{l}
a_{\th^i} \equiv a_{x^i \ot e^i} := 
\tilde{a}_{\th^i} \ot L_{e^i},\\
a_{\th^{\ast}_i} \equiv a_{x^{\ast}_i \ot e^{\ast}_i} :=
\tilde{a}_{\th^{\ast}_i} \ot L_{e^{\ast}_i},
\label{sup}
\ea
\ee
where $x^i, x^{\ast}_i \in s(\cw_{\G, c}(N))$ are given by (\ref{sge})
$e^i, e^{\ast}_i \in C_{\G, b}(N)$, 
$L : C_{\G, b}(N) \ra End(C_{\G, b}(N))$ is the left regular
representation of $C_{\G, b}(N)$ on itself, i.e. $L_{g} h = gh$ for
arbitrary $g, h \in C_{\G, b}(N)$. Obviously the above decomposition 
is not unique.
\section{Fock space for composite fermion and boson}
Let us describe the Fock space for composite fermion and
composite boson in details.
The the ground state vector $|0>$ is defined as usual, i.e. 
$a_{\th^{\ast}_i}|0> = 0$. We use here the notation
\be
a_{\th^i} \ |0> \equiv a^{+}_i |0> \equiv \th_i ,
\ee
and
\be
a^{+}_{\a} \ |0> 
\equiv (a^{+}_1)^{\a^1}\ldots(a^{+}_N)^{\a^N}|0>
\equiv (\th^{1})^{\a^1}\ldots (\th^N)^{\a^N}.
\ee
We use the similar notation for the operators $\tilde{a}_{x^i}$ 
and $L_{e^i}$. From relations (\ref{sup}) we obtain
\be
\ba{l}
a^{+}_{\a} \ |0> = \tilde{a}^{+}_{\a} \ot L_{\a} |0>\\
= (\th^{1})^{\a^1}\ldots (\th^N)^{\a^N}
\ot (e^{1})^{\a^1}\ldots (e^N)^{\a^N}.
\ea
\ee
Let us consider an example of reprepresentation corresponding
for an electron in singular magnetic field. We assume that the
electron is represented by one grassmann variable $\th$, i.e.
$\th^2 = 0$. For this representation we introduce here the 
following state vectors
\be
a^{+}_i |0> := 
\left(
\ba{c}
0\\
\vdots\\
\th\\
\vdots\\
0
\ea
\right),
\ee 
where $\th$ is on the $i$--th row. 
In this way an arbitrary state vector can be given in the following 
form
\be
a^{+}_{\a} \ |0> := \left(
\ba{c}
\th^{\a^1}\\
\vdots\\
\th^{\a^N}
\ea
\right),
\ee 
for $\a = (\a^1,\ldots,\a^N)$, $\a^i = 0$ or $1$ for $i = 1,\ldots,N$.
Now let us study this representation in more details. 
Observe that for even $N$ we obtain 
\be
\ba{cc}
\th_i \ \th_j = \th_j \ \th_i,&\;\th_i^2 = 0.
\ea
\ee
We have
\be
\th_i \th_j = \left(
\ba{c}
0\\
\vdots\\
\th\\
\vdots\\
\th\\
\vdots\\
0 
\ea
\right) = \th_j \th_i.
\ee
and
\be
\th_i^2 = \left(
\ba{c}
0\\
\vdots\\
\th^2\\
\vdots\\
0
\ea
\right) = 0.
\ee
Let us consider the case of $N = 2$
in more details. In this case we have the following states
\be
\ba{cc}
\th^1 = \left(
\ba{c}
\th\\
0
\ea
\right),&
\;  \th^2 = \left(
\ba{c}
0\\ 
\th
\ea
\right),
\label{sqa}
\ea
\ee
and
\be
\th^1 \ \th^2 = \left(
\ba{c}
\th\\ 
\th
\ea
\right).
\label{haf}
\ee
The filling factor for all these states (\ref{sqa}) and (\ref{haf})
is $v = \frac{1}{2}$. These two states (\ref{sqa}) contain quasiholes. 
Observe that the state (\ref{haf}) describe electron bond to two
fluxes which is transformed into a system containing two quasiparticles
and two fluxes. This system is nonlocal due to the anticommutativity
of fluxes. Note that the state (\ref{haf}) is the unique state 
corresponding to the filling factor $v = \frac{1}{2}$ without quasiholes. 
Observe that this state is completely filled by quasiparticles. This 
means thatfor $B = 2\Phi_0$ we obtain that $B_{eff} = 0$. In this 
way the state (\ref{haf}) corespond to the $\frac{1}{2}$--anomaly.

For odd $N$ we obtain
\be
\th_i \ \th_j = - \th_i \ \th_j \;\;\;\mbox{for}\; i \neq j .
\ee
For $i = j$ we obtain the identity $\th^i \th^i = \th^i \th^i$.
In this case the space of states can also be represented by 
the variable $\th$ such that we have
\be
\th^i = \left(
\ba{c}
0\\
\vdots\\
\th\\
\vdots\\
0
\ea
\right),
\ee
where $\th$ is on the $i$-th row. We have
\be
\th^i \th^j = \left(
\ba{c}
0\\
\vdots\\
\th\\
\vdots\\
\th\\
\vdots\\
0
\ea
\right)
= - \th^j \th^i.
\ee
This means that
\be
\th^i \ \th^j = 0
\ee
for all $i \neq j$. Observe that the quantum state 
$x^i x^j, \;(i \neq j)$ 
corresponding to particle coupled to magnetic fluxes 
at two different points disappear
\be
\th^i \ \th^j = x^i \ x^j \ot e^i \ e^j = 0,
\ee
and the state describing the particle coupled to a few fluxes at
the same point is also impossible. In fact we have
\be
(\th^i)^2 = (x_i)^2 \ot (e_i)^2 = (x_i)^2 \ot 1 = (x_i)^2.
\ee
This means that the state $(\th^i)^2$ is not equipped with a flux.
Let us consider as an example the case of $N = 3$ i.e. the filling
factor is $v = \frac{1}{3}$. In this case we have the states
\be
\ba{ccc}
\th^1 = \left(
\ba{c}
\th\\
0\\
0
\ea
\right),
&\;  \th^2 = \left(
\ba{c}
0\\ 
\th\\
0
\ea
\right),
&
\; \th^3 = \left(
\ba{c}
0\\
0\\
\th
\ea
\right),
\label{star}
\ea
\ee
which contain two quasiholes.
Observe that the following states
\be
\ba{ccc}
\th_1 \th_2 ,&\;  \th_1 \th_3,&\; \th_2 \th_3 .
\label{atar}
\ea
\ee
which contain one quasihole and the state
\be
\ba{c}
\th_1 \th_2 \th_3
\label{ful}
\ea
\ee
which not contain quasiholes are impossible. Hence in this case the 
single quasiparticle states with two quasiholes are possible! 
This also means that the Landau lowest level can not be completely
filled by quasiparticles but fractionally filled only! 
This is just the FQHE.

We can see that there is a well defined scalar product on $\ca$ 
\cite{ral}
\be
<\th^{\sigma}|\th^{\tau}>_{q,b} = \Sigma_{\pi \in S_n} \
 \chi_n (\pi) \ <\th^{i_1}|\th^{\pi(j_1)}> \ldots 
<\th^{i_n}|\th^{\pi(j_n)}>,
\ee
where $\th^{\sigma}, \th^{\tau} \in \ca^n$,
$\th^{\sigma} = \th^{i_1} \ldots \th^{i_n}$,
$\th^{\tau} = \th^{j_1} \ldots \th^{j_n}$, and
\be
\chi (\pi) := \Pi_{(i,j)\in J(\pi)} \ c_{ij}
= \Pi_{(i,j)\in J(\pi)} \ s_i^{n_i} \ b_{ij} ,
\ee
$J(\pi) := \{(i,j): 1 \underline{<} i \underline{<} j \underline{<} n\},
\pi (i) > \pi (j)$, $n_i$ is the number of elements of the set
$K(\pi) := \{(i,j) \in J(\pi): i = j\}$.
It follows from the theorem 
of Bo$\dot{z}$ejko and Speicher \cite{bs2}
that the corresponding scalar product is positive definite.
\section{Appendix}
Let us shortly recall the concept of commutation factors on
an arbitrary abelian group $\G$, \cite{sch}.
A mapping $c: \G \times \G \ra I\!\!\!\! {\bf C}$ such that
\be
\ba{l}
c(\a+\b,\g) = c(\a,\g)c(\a,\g),\;\;\;  c(\a,\b+\g) = c(\a,\b)c(\b\g),\\
\ea
\ee
is said to be a bicharacter on $\G$. If in addition we have the relation
\be
c(\a,\b)c(\b,\a)=1,
\ee
then $c$ is said to be a commutation factor on $\G$.
We use here the following notation
$c(\si_i, \si_i) \equiv c_{ii} \equiv q_i$ and 
$c(\si_i, \si_j) = c_{ij}$ for $i \neq j$,
where $q_i = +1$ or $-1$ is said to be a parity of $c$, 
$c_{ij} \in I\!\!\!\! {\bf C} \setminus \{0\}$ are deformation
parameters such that $c_{ij} c_{ji} = 1$;
$\a = \Sigma_i \a^i \si_i$, $\b = \Sigma_j \a^j \si_j$, 
$\si_i$ are generators of $\G$. 
Note that we have in general the following formula for the
commutation factor $c$ on the group $\G$ 
\be
\ba{c}
c(\a, \b) = \Pi_{i,j} (c_{ij})^{\a^i \b^j}\\
= \Pi_{i} (c_{ii})^{\a^i\b^i} \Pi_{i<j} (c_{ij})^{\a^i\b^j - \a^j\b^i}.
\label{don}
\ea
\ee
It is obvious that the set $\G_0 := \{\a \in \G : c(\a, \a) = 1\}$
is a subgroup of $\G$ of index at most $2$. This means that the
quotient $\G/\G_0$ is isomorphic to the group $Z_2$ or is trivial.
The commutation factor $c$ can be given in the following form
\be
c(\a, \b) = c'(\a, \b) \ b(\a, \b),
\ee
where $b$ and $c'$ are two new commutation factors on $\G$. 
The factor $b$ is defined by the relations
\be
b(\si_i, \si_j) \equiv b_{ij} := \left\{
\ba{rl}
+ 1&\;\mbox{for};\ i = j\\
-c_{ij}&\;\mbox{for}\; i \neq j \; \mbox{if} \; q_i = q_j = -1 \\
 c_{ij}&\;\mbox{in the remaining cases}
\ea
\right..
\ee
It follows immediately from the formula (\ref{don}) that
\be
\ba{l}
b(\a, \b) = \Pi_{i<j} (b_{ij})^{\a^i\b^j - \a^j\b^i}.
\label{ncf}
\ea
\ee
The factor $c'$ is defined by the formula
$c'_{ii} \equiv c_{ii} \equiv q_i$ and 
\be
c'(\si_i, \si_j) \equiv c'_{ij} = \left\{
\ba{ll}
- 1&\mbox{if}\;q_i = q_j = -1, \\
+ 1&\mbox{otherwise}
\ea
\right.
\label{sad}
\ee
for $i \neq j$. It is interesting that the factor $c'$ can be reduced
to the group $\pi(\G) \equiv \G/\G_0$, where $\pi : \G \ra \G/\G_0$ 
is the quotient map. We can see that there is the relation
\be
c'(\a, \b) = (-1)^{\pi(\a)\pi(\b)}.
\ee
For $\G = Z^N \equiv Z \oplus \ldots Z$ ($N$--sumands) we have
\be
\ba{ll}
q_i = (-1)^{S_i},& \; b_{ij} = \omega^{\Omega_{ij}}, 
\label{omo}
\ea
\ee
where $S_i = 0, 1$, and $\Omega_{ij}$ are elements of a skew-symmetric 
integer-valued, matrix, and $\omega \neq -1$ is some complex parameter 
\cite{zoz}. If $\omega = exp(\frac{2 \pi i}{n})$, $n \underline{>} 3$ 
and $S_{i} \equiv 0$, then the grading group $G \equiv Z^N$ can 
be reduced to $\G = Z^N_n \equiv Z_n \oplus...\oplus Z_n$. 
If $\om = -1$, then the grading group $\G$ can be reduced to 
the group 
\be
\ba{c}
G \equiv Z_2^N := Z_2 \oplus ... \oplus Z_2\;\;\;(N sumands).
\label{gru}
\ea
\ee
In follows from our considerations in the Section 2 that for the
charged particle in the singular magnetic field we have
$S_{i} = N+1$ (mod $2$), and
\be
\Om_{ij} :=
\left\{
\ba{lcc}
0&\mbox{for}&i=j,\\
1&\mbox{for}&i\neq j .
\ea
\right.
\ee
This means that 
\be
b_{ij} :=
\left\{
\ba{lcc}
+1&\mbox{for}&i=j\\
-1&\mbox{for}&i\neq j 
\ea
\right. ,
\ee
and the factor $c$ is given by the formula
\be
\ba{c}
c_{ij} = -(-1)^{N} (-1)^{\Om_{ij}}.
\label{com}
\ea
\ee
In such case if the number $N$ of fluxes is even (i. e. for composite 
fermions) we obtain
\be
c_{ij} = -(-1)^{\Om_{ij}} =
\left\{
\ba{lcc}
-1&\mbox{for}&i=j\\
+1&\mbox{for}&i\neq j 
\ea
\right. .
\label{fom}
\ee
If $N$ is odd (composite bosons), then we obtain 
\be
\ba{c}
c_{ij} = (-1)^{\Om_{ij}} :=
\left\{
\ba{lcc}
+1&\mbox{for}&i=j,\\
-1&\mbox{for}&i\neq j .
\ea
\right.
\label{wom}
\ea
\ee

\vspace{1cm}

{\bf Acknowledgments}\\
The author would like to thank to J. Lukierski for discussion, 
to R. Gielerak for some remarks, and to A. Borowiec for any others 
help.

\end{document}